\documentclass[12pt,twoside,reqno]{amsart}
\usepackage[colorlinks=true,citecolor=blue]{hyperref}
\usepackage{mathptmx, amsmath, amssymb, amsfonts, amsthm, mathptmx, enumerate, color}
\usepackage{graphicx}
\usepackage{epstopdf}
\usepackage{xcolor}

\theoremstyle{definition}

\numberwithin{equation}{section}

\begin{document}
\setcounter{page}{1}
\vspace*{1.0cm}
\title[Avoiding Traps]
{Avoiding Traps in Nonconvex Problems}
\author[S. Deyo, V. Elser]{ Sean Deyo$^{1,*}$, Veit Elser$^1$}
\maketitle
\vspace*{-0.6cm}

\begin{center}
{\footnotesize {\it

$^1$Department of Physics, Cornell University, Ithaca, NY, USA

}}\end{center}

\vskip 4mm {\small\noindent {\bf Abstract.}
Iterative projection methods may become trapped at non-solutions when the constraint sets are nonconvex. Two kinds of parameters are available to help avoid this behavior and this study gives examples of both. The first kind of parameter, called a hyperparameter, includes any kind of parameter that appears in the definition of the iteration rule itself. The second kind comprises metric parameters in the definition of the constraint sets, a feature that arises when the problem to be solved has two or more kinds of variables. Through examples we show the importance of properly tuning both kinds of parameters and offer heuristic interpretations of the observed behavior.

\noindent {\bf Keywords.}
Projection methods; Fixed-point algorithms; Nonconvex problems; Logical satisfiability; Dominating sets; Machine learning. }

\renewcommand{\thefootnote}{}
\footnotetext{ $^*$Corresponding author.
\par
E-mail addresses: sjd257@cornell.edu (S. Deyo), ve10@cornell.edu (V. Elser).
\par
Submitted May 10, 2021. Revised November 4, 2021.}

\section{Introduction}

Iterative algorithms whose elementary operations are projections to constraint sets perform well on many nonconvex problems for which one does not have the convergence guarantees one has in convex problems. For some combinatorially hard problems these algorithms routinely outperform state-of-the-art algorithms that find solutions by exhaustive search \cite{bitretrieval}. This success, or the apparent ability of the projection-based search to very significantly reduce the size of the space being searched, is poorly understood.

The standard criteria for evaluating iterative projection algorithms in the convex case do not apply in the nonconvex case. The amount of time the algorithm spends refining a solution, once its locally convex basin has been encountered, is negligible compared to the time needed to find the basin. Efficiency in basin discovery completely overshadows the benefits of good convergence within the basin.

Iterative projection algorithms usually follow a fixed-point principle, where fixed-points imply a solution, but these algorithms may still get trapped on non-solutions in a dynamic sense \cite{global}. When this happens, the iterations meander indefinitely within a small domain far from a true fixed point. Eliminating or mitigating this behavior by tuning the parameters of the algorithm is the focus of this paper.

Most iterative projection algorithms have parameters that apply to the general case, independent of the application's constraint sets. For such parameters we use the machine learning term ``hyperparameter." Relaxation parameters are examples of hyperparameters. Another type of parameter has only come to light more recently, in applications that require multiple kinds of variables \cite{matrixproduct}. These applications have variable-scaling freedom that is not a geometrical isometry and therefore changes the algorithm through its effect on the projections. We refer to such parameters as ``metric parameters."

\section{Hyperparameters and metric parameters}\label{sec:hypervsmetric}

One of the best known hyperparameters is the relaxation that is often applied to the standard Douglas-Rachford iteration
\begin{equation}\label{RRR}
x\mapsto (1-\beta/2)\,x+(\beta/2)\,R_B(R_A(x)),
\end{equation}
with $\beta \in\; ]0,2[$. This is called the generalized Douglas-Rachford method~\cite{tutorial}. Here $R_A$ and $R_B$ are the reflectors for the constraint sets $A$ and $B$: $R_i(x)=2P_i(x)-x$, where $P_i(x)$ is the projection of $x$ to set $i$. In the case where one of the sets $A$ or $B$ is nonconvex and there exist points $x$ where the projection is not unique, we use the fact that such points have measure zero \cite{nonuniqueProj} and are not expected to arise in computations with sufficiently high floating point precision.
The hyperparameter $\beta$ is the same parameter that independently was deemed important when this iteration was proposed --- by an engineer unaware of Douglas-Rachford --- for the phase retrieval application \cite{HIO}. To make the point that the scope of hyperparameters is broad, we give an example in section \ref{sec:doublereflect} of a very different generalization of Douglas-Rachford, in which tuning a hyperparameter is critical for success.

In imaging applications \cite{benchmarks}, or even sudoku \cite{sudoku}, where there is just one kind of pixel or cell, the question of metric parameters never came up. An application where introducing metric parameters makes sense is non-negative matrix factorization (NMF) \cite{NMF}.
In NMF one seeks a low-rank factorization of a rectangular non-negative matrix $Z=X Y$, where the factors are themselves non-negative. The rows of $Z$ may be interpreted as data vectors that can be expressed as non-negative mixtures, given by $X$, of a set of non-negative features, the rows of $Y$. NMF is non-unique with respect to rescaling ($X\to XD$, $Y\to D^{-1}Y$, arbitrary diagonal and positive $D$). To remove this ambiguity and also to make the problem compact, we can choose to impose a norm on the columns of $X$ or the rows of $Y$. If we decide to normalize the features ($Y$), what setting of the norm do we choose and why is that a choice of metric?

To motivate the term ``metric," consider the standard distance that would be used in defining the projections for the NMF application:
\begin{equation}\label{eq:metric}
d\left((X,Y),(X',Y')\right)=\sqrt{\|X-X'\|_2^2+\|Y-Y'\|_2^2}.
\end{equation}
Now, if we chose to impose normalization on each row vector $y$ of $Y$, say $\|y\|_2=\eta$, we could alternatively work with the rescaled variables $\tilde{Y}=Y/\eta$, normalization constraint $\|\tilde{y}\|_2=1$, and the distance
\begin{equation}\label{eq:scaledmetric}
d\left((X,\tilde{Y}),(X',\tilde{Y}')\right)=\sqrt{\|X-X'\|_2^2+\eta^2 \|\tilde{Y}-\tilde{Y}'\|_2^2}.
\end{equation}
This rewriting by a parameterized distance provides an interpretation of the norm parameter $\eta$. Consider the projection to the bilinear constraint, $X\tilde{Y}=Z/\eta$. From the distance \eqref{eq:scaledmetric} we see that for small $\eta$ the $\tilde{Y}$ variables (features) are more compliant than $X$ (mixtures), and we should expect that product-constraint inconsistencies are resolved mostly by changing the features. The opposite, or more reliance on changing the mixtures, is expected for large $\eta$. Whether one chooses to normalize $y$ to $\eta$ and use the standard metric~\eqref{eq:metric} or keep $y$ normalized to $1$ and use the modified metric~\eqref{eq:scaledmetric} is a matter of implementation. In this paper our equations assume the former choice.

NMF is one of the earliest techniques of machine learning, and we anticipate that variable-type metric sensitivity will grow in relevance as projection methods find their way into this domain$^1$\footnote{1. The term ``parameter" is potentially confusing in the machine learning context. For example, the \textit{weight para\-meters} of a neural network that are learned from data are variables from the perspective of the optimization algorithm. The latter may use \textit{metric parameters} to more efficiently optimize the \textit{weight variables}.}. In particular, splitting methods lend themselves naturally to optimization on networks, where variables appear on both nodes and edges of the network and are clearly dissimilar.

\section{Hyperparameters}

\subsection*{A tuneable double-reflector algorithm}\label{sec:doublereflect}

In this section we consider a generalization of the Douglas-Rachford iteration that is very different from the standard relaxation \eqref{RRR} :
\begin{equation}\label{doubleref}
x\mapsto \mathrm{DR}_n[\delta](x)=\frac{1}{1+n}\;\sum_{r=0}^{n}\;(R_B[\delta] \circ R_A[\delta])^r(x)\;.
\end{equation}
The number of double reflections, $n$, is one of the algorithm's hyperparameters. However, we will see that this algorithm only succeeds when the reflectors are themselves parameterized,

\begin{equation}
R_i[\delta](x)=(2-\delta)\,P_i(x)-(1-\delta)\,x\;,
\end{equation}
with $\delta\in\;[0,1]$. The original Douglas-Rachford iteration is recovered for $n=1$ and $\delta=0$ and will be referred to as $\mathrm{DR}_1[0]$.

The idea behind \eqref{doubleref} is that the double reflector acts as the identity when $x$ is near a feasible point, and therefore the average of any number of double reflections fixes $x$. On the other hand, multiple ($n>1$) applications of the double reflector might be better at ejecting $x$ from a trap when it is not near a feasible point. Through elementary analysis and numerical experiments we will argue that there is an optimal $\delta^*$ such that when $\delta>\delta^*$, and the reflections are in effect contractive, trapping the iterations at non-solutions, while for $\delta<\delta^*$ the iterate diffuses too freely to notice even the true solutions. The algorithm works best when $\delta$ is tuned to the dynamical transition point $\delta^*$.

Hard feasibility problems, including the one we feature below, can often be formulated where one constraint set, say $A$, is finite and the other, $B$, is a hyperplane. Traps arise when a point $a\in A$ is very close to $B$, that is, when the distance $\Delta=\|a-b\|$, to the proximal point $b=P_B(a)$ on the hyperplane, is very small. 
We analyze the local trapping/escaping behavior by replacing the sets $A$ and $B$ by proximal points $a\in A$,  $b\in B$ and the ambient space by the line passing through these points. Let $x$ be the (1D) coordinate along this line with $a$ corresponding to $x=0$ and $b$ corresponding to $x=-\Delta$. In this simplified model of the local behavior one finds
\begin{equation}\label{doubleref1D}
\mathrm{DR}_n[\delta](x)=\gamma\, x+c\;,
\end{equation}
where
\begin{align}
\gamma&=1-\frac{n(1-q)-q(1-q^{n})}{(1+n)(1-q)}\nonumber\\
q&=(1-\delta)^2\nonumber\\
c&=(1-\gamma)\left(\frac{1-\delta}{\delta}\right)\Delta\; .\nonumber
\end{align}
For comparison,
\begin{equation}
\mathrm{DR}_1[0](x)=x+\Delta
\end{equation}
represents a step-wise escape, where $O(1/\Delta)$ iterations are needed before $x$ has changed by $O(1)$. This is an estimate of the number of iterations, in the original problem, for $P_A(x)$ to be significantly different from the point $a$ of the trap.

Since $0<\gamma<1$ for $0<\delta<1$, iteration \eqref{doubleref1D} is contractive and looks problematic because there is always a (non-solution) fixed point:
\begin{equation}\label{1Dfixedpoint}
x^*=\frac{c}{1-\gamma}=\left(\frac{1-\delta}{\delta}\right)\Delta\;.
\end{equation}
However, \eqref{1Dfixedpoint} should be seen as instructions on the proper use of the algorithm. Since one constraint set is finite there will be a smallest $\Delta$ that poses the greatest trapping risk. But by setting $\delta=\delta^*=O(\Delta)$, the fixed point $x^*$ of the 1D dynamics is sufficiently far from the trap that $P_A(x^*)$ will likely be different from the original trapping point $a$. Since $c\sim n\,\Delta$ for $\delta\to 0$, a single iteration of $\mathrm{DR}_n$ is roughly the same as $n$ iterations of $\mathrm{DR}_1$, although both schemes require $n$ double-reflector computations.

\subsection*{Experiments with logical satisfiability}

To illustrate the effect of the hyperparameter $\delta$ in a setting with potentially many traps, we turn to the logical satisfiability problem (SAT). In SAT we have a set $C$ of \textit{clauses} and a set $V$ of \textit{variables}. Thinking of these as vertices of a bipartite graph $G$, the search variables $E$ in our constraint formulation correspond to edges $c\to v$, where $c\in C$ and $v\in V$. In the SAT interpretation, the variable-vertices $v$ incident on a particular $c$ in $G$ correspond to the Boolean variables that participate in one clause of a logical formula in conjunctive normal form. The clause itself is a disjunction
\begin{equation}
y_c=\bigvee_{v\in V:\; c\to v\;\in\; E}n_{c\to v}\circ x_{c\to v}
\end{equation}
where the $x$ are Boolean variables, $n\circ$ specifies whether to apply negation, and $y_c$ is the Boolean value of clause $c$. The object in SAT is to find an assignment to the $x$ such that the conjunction
\begin{equation}\label{conj}
\bigwedge_{c\in C}y_c
\end{equation}
is true. The set of Boolean variables $x_{c\to v}$ incident on the same $v$ (but appearing in different clauses $c$) should all be equal. However, in the two-constraint formulation \cite{divideconcur} to which we turn next, these are treated as independent in one of the constraints.

The two constraint sets live in a space of dimension $|E|$. Set $A$ is finite and imposes the truth of each clause (otherwise the conjunction \eqref{conj} is false). We encode \textsc{True} and \textsc{False} for the Boolean variables as respectively $x=+1$ and $x=-1$ and use multiplication by $n=-1$ for negation:
\begin{equation}
\begin{array}{lll}
A:& &\\
& \forall\; c\to v\;\in E:& x_{c\to v}\in \{-1,1\}\;,\nonumber\\[10pt]
& \forall\; c\in C:& +1\in {\displaystyle \bigcup_{v\in V:\; c\to v\;\in\; E}}n_{c\to v}\;x_{c\to v}\;.\nonumber
\end{array}
\end{equation}
Because the $A$ constraint imposes the discreteness of the Boolean variables we are free to use the following relaxed, continuous constraint for $B$:
\begin{equation}
\begin{array}{lll}
B:& \forall\; v\in V:\; \exists\; \bar{x}_v\in \mathbb{R}\;:&\forall\; (c\in C : {c\to v}\in E): x_{c\to v}=\bar{x}_v\;.
\end{array}
\end{equation}

We consider a hard instance of 3-SAT, where each clause involves exactly three variables and there are altogether $|V|=500$ Boolean variables in the logical formula. The number of clauses $|C|=2100$ was tuned so that a typical random instance has roughly even odds of being satisfiable \cite{SATtransition} (our instance is satisfiable). All our results will be for iteration $\mathrm{DR}_3$, for which the number of double reflections ($n=3$) is large enough that a transition in behavior with $\delta$ is easy to discern. Each trial starts with a random initial point $x$, with each of the $|E|$ variables chosen in the range $[0,1]$.

\begin{table}[t]
\begin{tabular}{ c|r|c } 
 $\delta$ & successes/trials & iterations/solution\\
 \hline
 $0.001$ & $0/100$ & --- \\ 
 $0.002$ & $0/100$ & --- \\ 
 $0.005$ & $90/100$ & $3.23\times10^3$ \\ 
 $0.010$ & $86/100$ & $4.57\times10^3$ \\  
 $0.020$ & $35/100$ & $2.27\times10^4$ \\ 
 $0.050$ & $0/100$ & --- \\ 
 \hline
\end{tabular}
\vspace{.1in}
\caption{Performance statistics for several values of $\delta$ as algorithm \eqref{doubleref} tries to solve an instance of the 3-SAT problem with $500$ Boolean variables and $2100$ clauses. Each trial was capped at $10^4$ iterations.}
\label{tab:cnf}
\end{table}

\begin{figure}
  \centering
  \includegraphics[width=.44\textwidth]{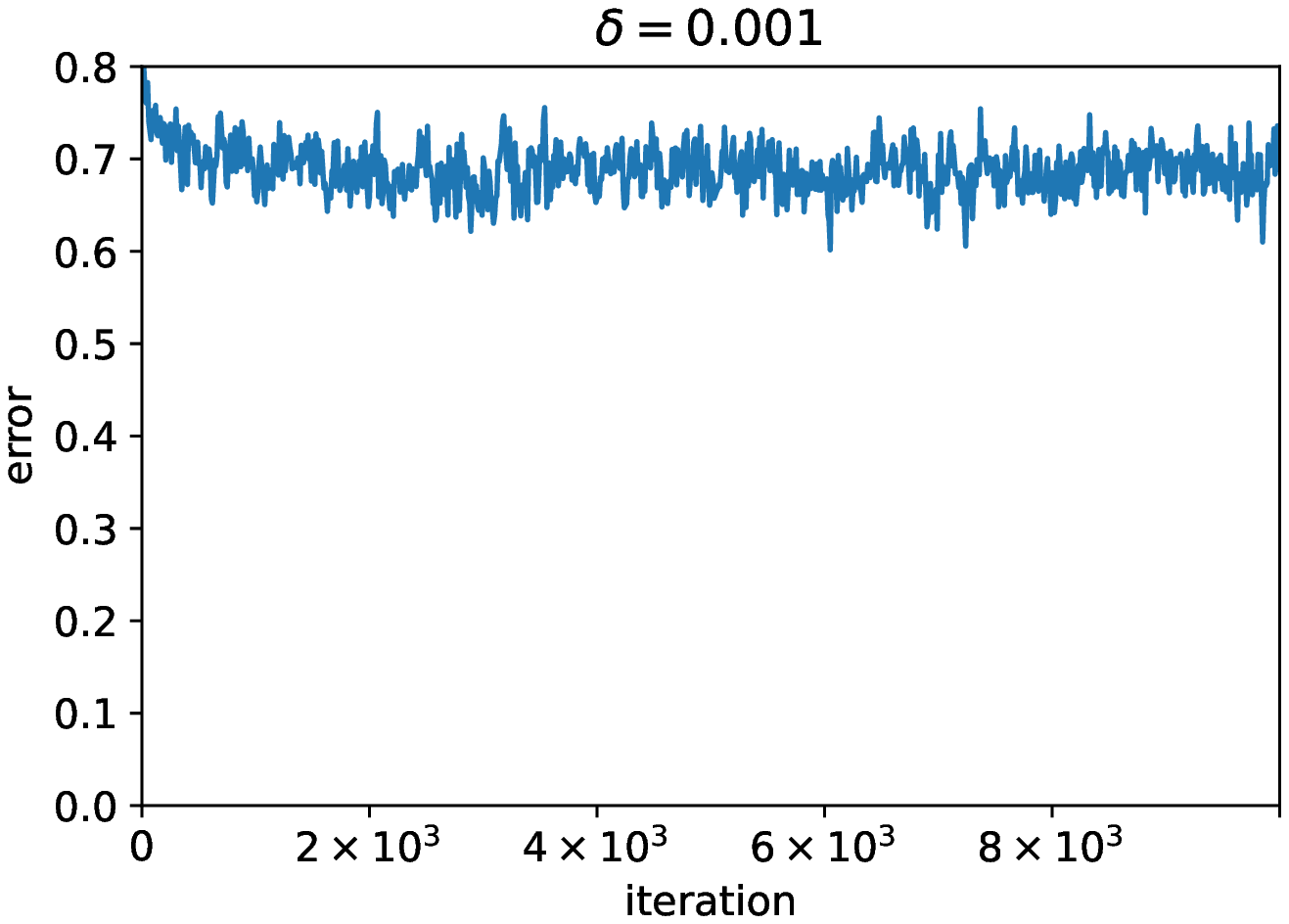}\quad
  \includegraphics[width=.44\textwidth]{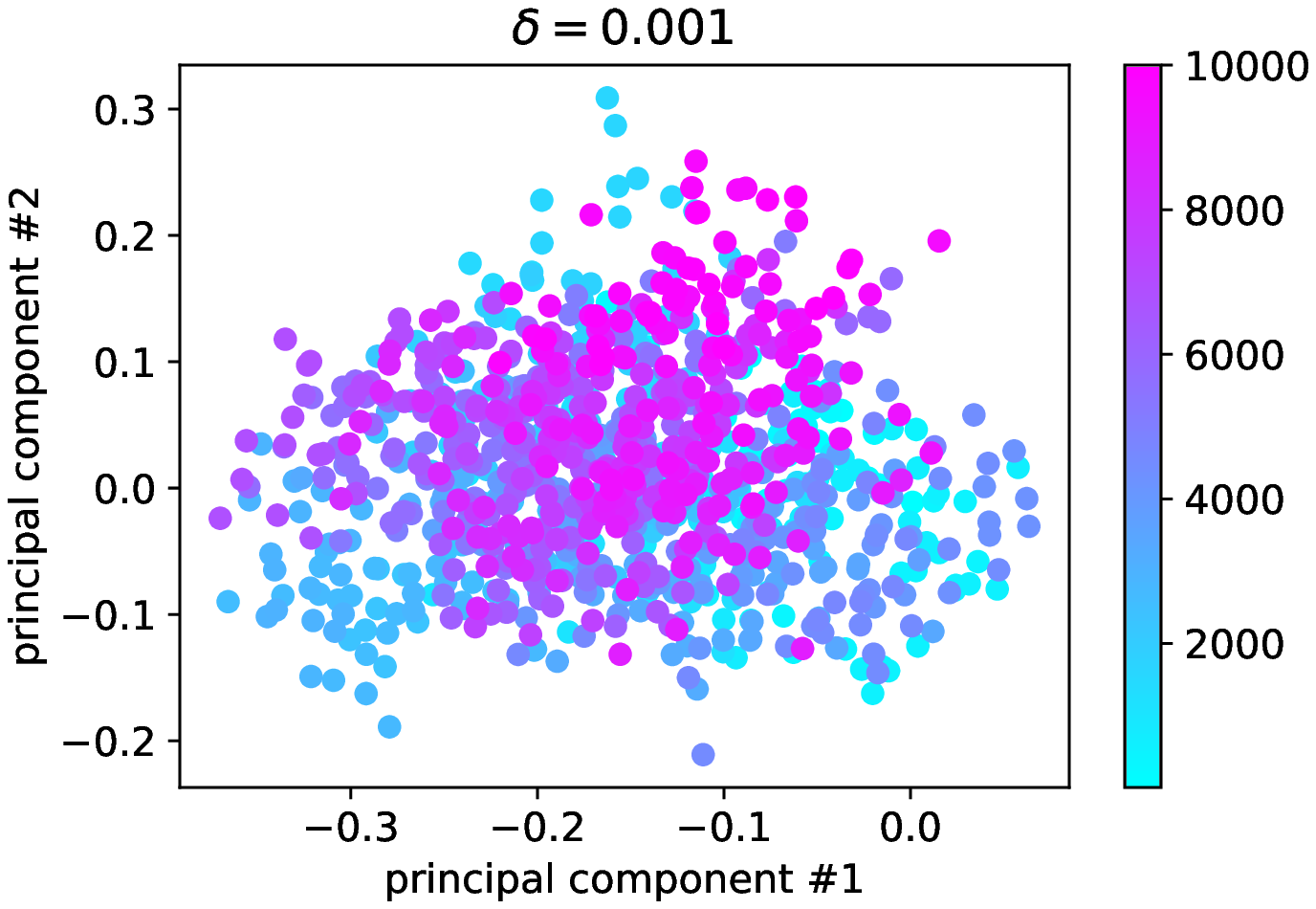}\\
  \includegraphics[width=.44\textwidth]{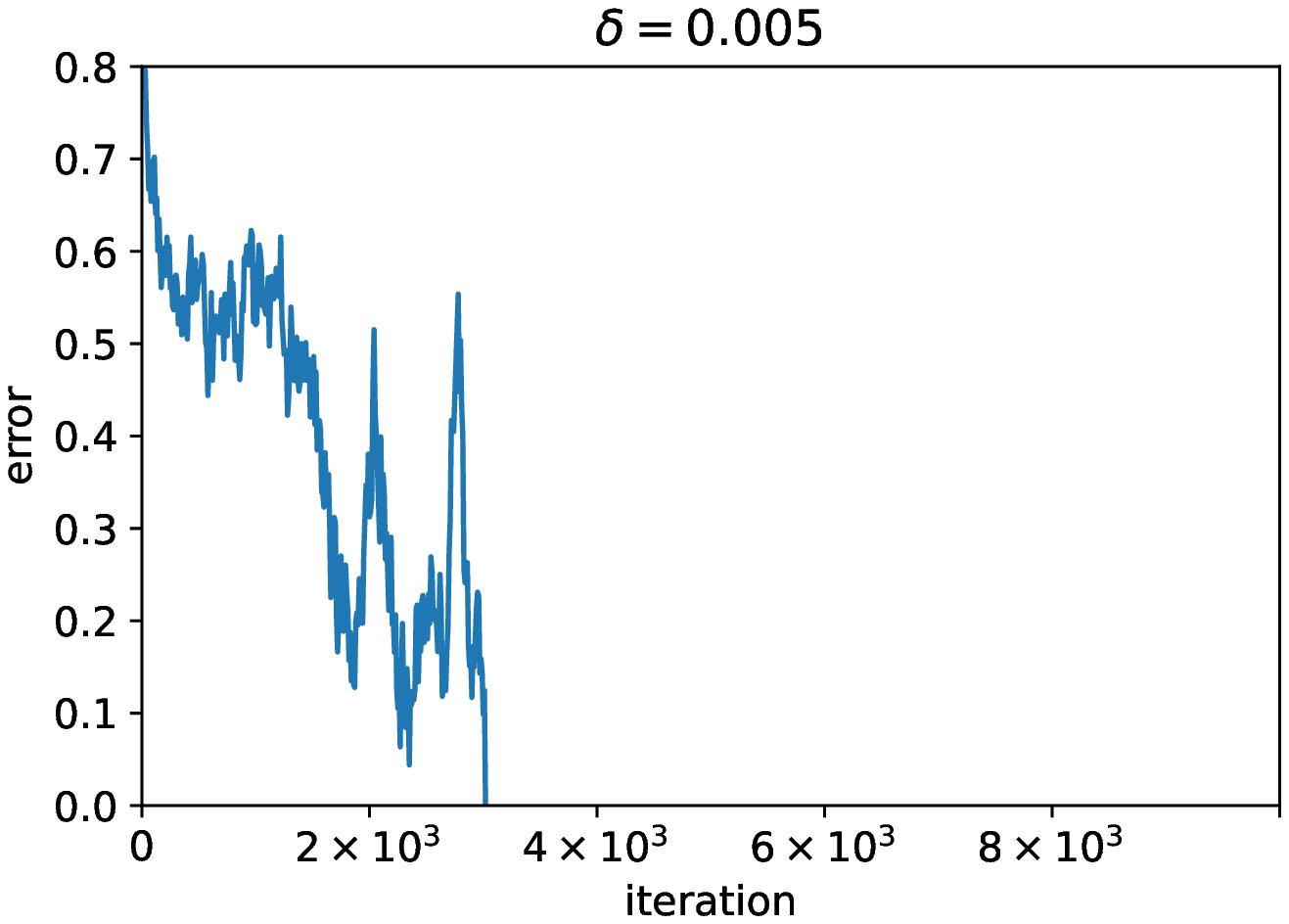}\quad
  \includegraphics[width=.44\textwidth]{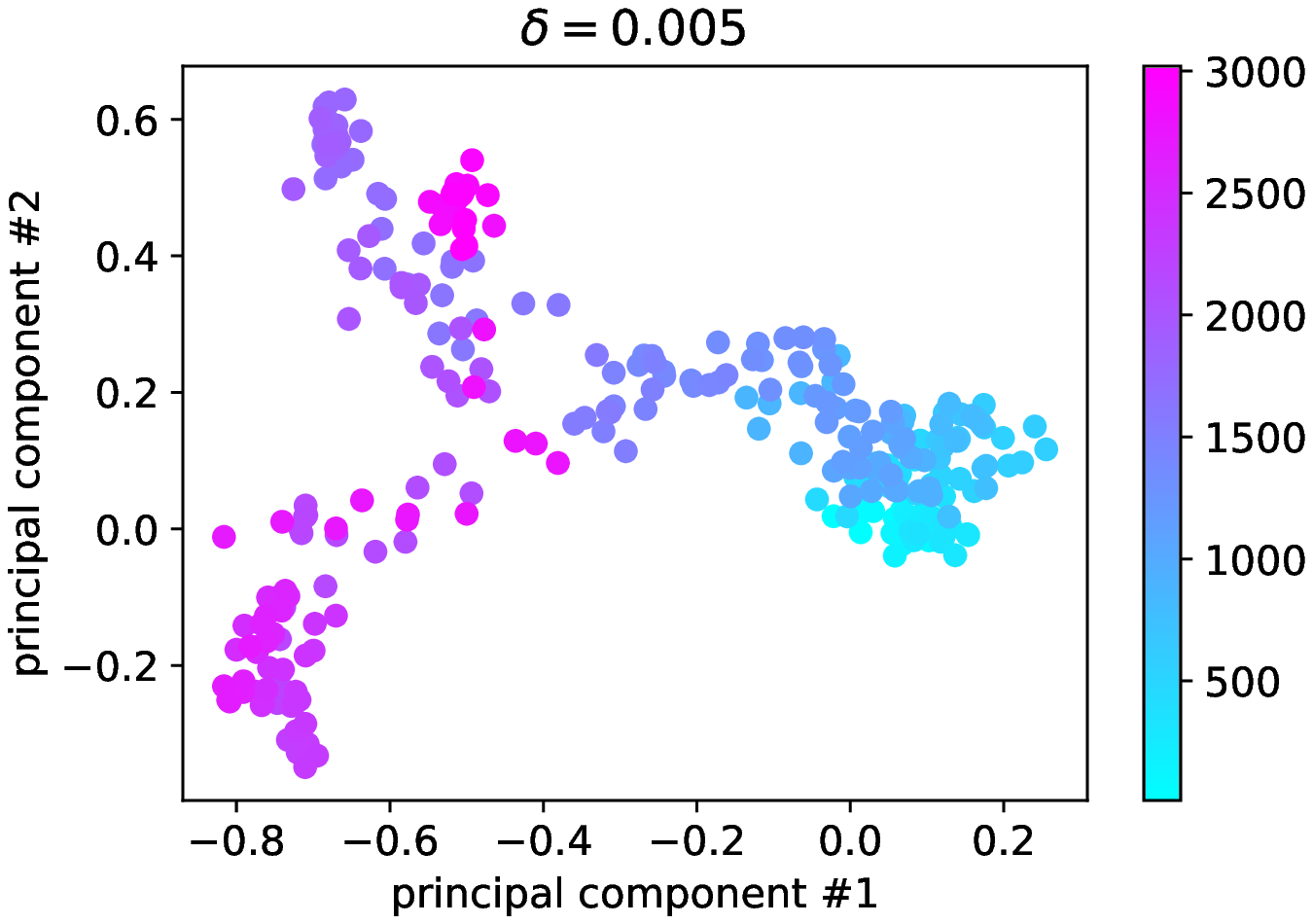}\\
  \includegraphics[width=.44\textwidth]{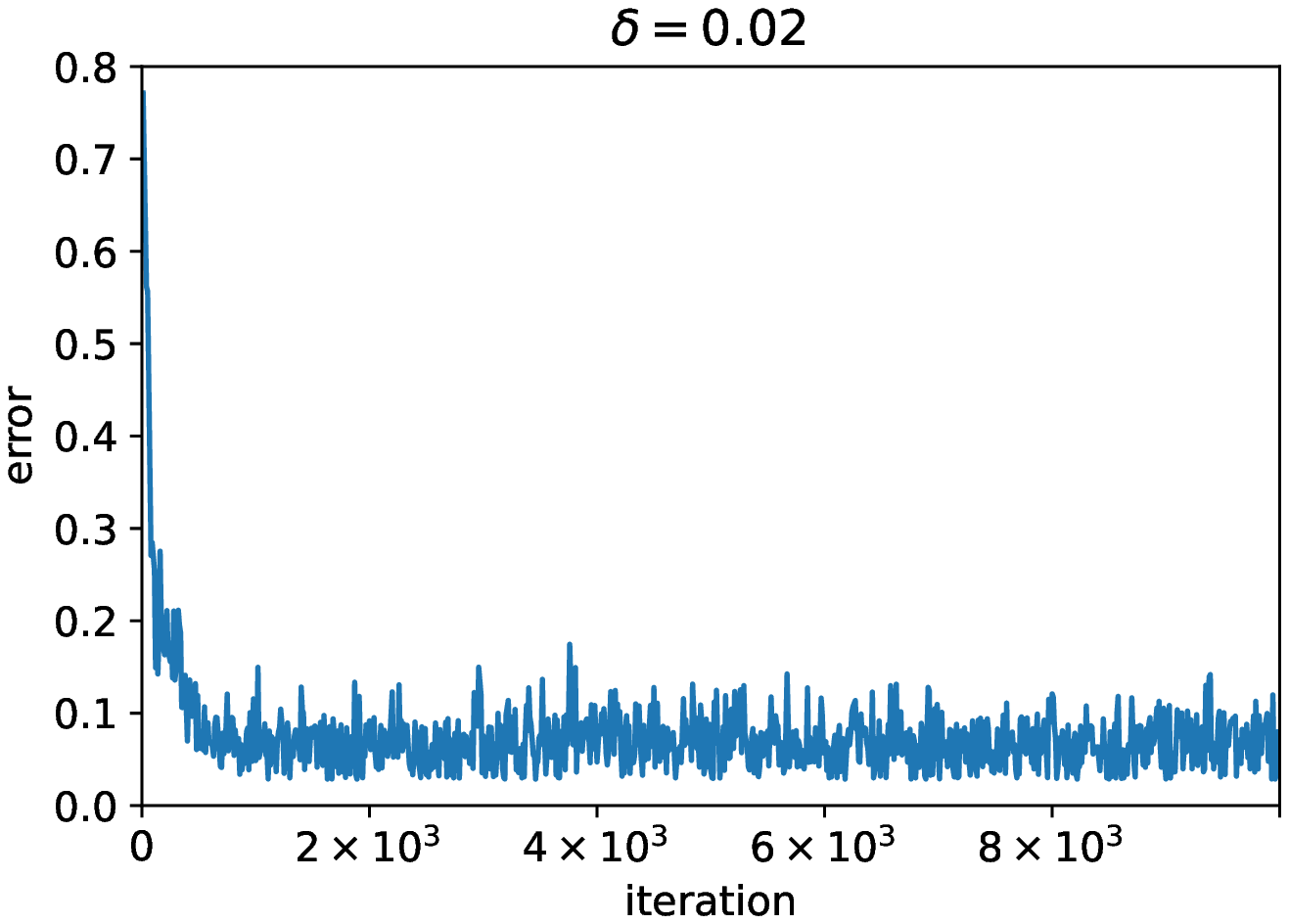}\quad
  \includegraphics[width=.44\textwidth]{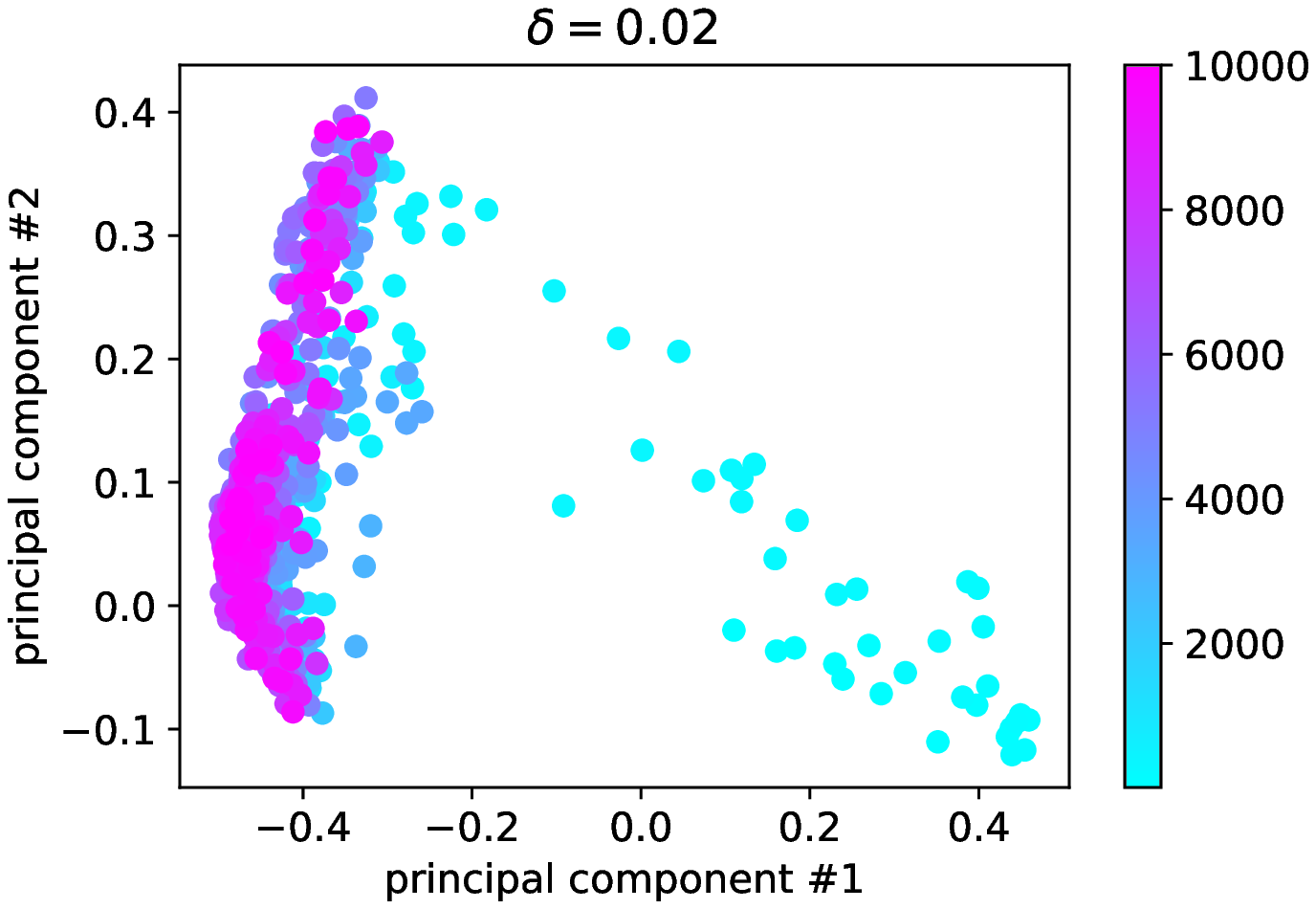}
  \caption{Plots of the constraint error $\epsilon$ (left column) and the iterate $x$ (right column) of the $\mathrm{DR}_3$ algorithm applied to an instance of 3-SAT. The plots on the right show the iterate $x$ projected into a 2D space by PCA. These data are plotted in increments of $10$ iterations up to a maximum of $10^4$ iterations, with the color bar indicating the iteration number.}
  \label{fig:cnf}
\end{figure}

The performance statistics are given in Table \ref{tab:cnf}. The smallest values of $\delta$ are completely ineffective. There is a sharp onset of good performance at $\delta=0.005$, but the larger values of $\delta$ are also ineffective. To understand why, it is helpful to look more closely at a few individual trials. To visualize the behavior, we store the time series of the $|E|$ search variables in a matrix, where each row specifies a point in $\mathbb{R}^{|E|}$. We then find the principal component axes for this matrix and project each row onto the first two principal components to obtain the 2-dimensional plots in Figure \ref{fig:cnf} of the search as a function of time. Since we are projecting points with root-mean-square distance $O(\sqrt{|E|})$, we also divide the principal components by $\sqrt{|E|}$ to normalize the length scale in the 2D plots. Alongside each PCA time series we provide the time series of the constraint error 
\begin{equation}
\epsilon=|P_A(x)-P_B(R_A(x))|
\end{equation}
for the points $x$ generated by the $\mathrm{DR}_3$ iteration.

In Figure \ref{fig:cnf} we see that the search with  $\delta=0.001$ (top plots) jumps around aimlessly with large $\epsilon$ and never finds a solution. Our interpretation is that when $\delta$ is small, each reflection $R_i(x)$ is almost the same as the pure reflection $2P_i(x)-x$ and the algorithm never allows itself to fall into a basin of lower $\epsilon$. When $\delta=0.02$ (bottom plots of Fig. \ref{fig:cnf}), the algorithm quickly gets trapped on the first basin it finds, corresponding to a non-solution, and remains there for the rest of the search. The intermediate choice $\delta=0.005$ (middle plots of Fig. \ref{fig:cnf}) seems to be just right. The algorithm now appears to start exploring a basin once or twice but only stays for a few hundred iterations before wandering to another basin. After about $3\times10^3$ iterations it finds a basin that has a solution and converges toward it.

\section{Metric parameters}

Metric parameters, unlike hyperparameters, are special for each intended application and should always be considered when there is more than one type of variable. Below we give two examples, both involving applications where the variables take discrete values and live on a network. As with hyperparameters, the settings of the metric parameters make all the difference between an algorithm that never finds solutions and one that does so consistently.

Often a metric parameter will not have an obvious interpretation, or will have non-obvious interactions with the other metric parameters. Tuning these metric parameters can be done by hand and is informed by appropriate diagnostics that go beyond the standard ``success rate" statistic. However, an automated procedure to expedite this process is desirable, especially when the number of metric parameters is large. We use a scheme where the current status of the search informs the rule for the parameter updates, and apply these updates adiabatically so as not to upend the fixed-point properties of the algorithm being used.

Our metric parameter update rule is based on the following heuristic. We want to prevent the algorithm from getting stuck on a partial solution, wherein some variables hardly change between the $A$ and $B$ projections while others are changing very much. We avoid this by giving a smaller metric parameter to the variables that are hardly changing, thereby lowering the penalty for changing them when the next projection comes, and vice versa for variables that change too much.

To be precise, suppose we partition the variables into $k$ types, with metric parameters $\eta_1$, $\eta_2$, ... $\eta_k$. Let $l_i$ denote the number of variables of type $i$, so that the complete variable vector can be written as a concatenation $x=(x_1,x_2,...,x_k)$, where $x_1$ is a vector of length $l_1$ and so on. The projections are made using the distance function 
\begin{equation}
    d(x,x')=\sqrt{\sum_{i=1}^k \|x_i-x'_i\|^2}.
\end{equation}

In each iteration, we compute the normalized rms error 
\begin{equation}
    \epsilon_i=\frac{1}{\eta_i}\sqrt{\|P_A(x)_i - P_B(R_A(x))_i\|^2/l_i}
\end{equation}
for each variable type. We then compare each $\epsilon_i$ to the average $\epsilon = \sqrt{\frac1k \sum_i \epsilon_i^2}$ and adjust the metric parameters according to
\begin{equation}
\eta_i \to \eta_i \left(1+\alpha\left(\epsilon_i/\epsilon-1\right)\right)
\label{update}
\end{equation}
where $\alpha\ll1$ is a small but positive tuning parameter. Alternatively, since only the relative weights of the variable types matter, one can set $\eta_1=1$ and take
\begin{equation}
\eta_i \to \eta_i \left(1+\alpha\left(\epsilon_i/\epsilon_1-1\right)\right)
\label{update1}
\end{equation}
for $i>1$. 

Automatically updating the metric parameters in this fashion encourages the smaller $\epsilon_i$'s to become larger and vice versa, which leads to the $\epsilon_i$'s being highly correlated, which in turn generally leads to more successful searching. This approach also replaces the problem of tuning some (possibly large) number of metric parameters with the simpler question of choosing a value of $\alpha$. One must have $\alpha\ll1$ in order to make the metric parameter updates adiabatic and thereby preserve the local convergence properties of the algorithm. One also wants $\alpha\gg1/N$, where $N$ is the total number of iterations to be run, so as to accomplish the desired tuning within the intended length of the run. Accordingly, for large $N$ there can be a rather generous range of $\alpha$ that will work.

\subsection*{Experiments}
For the following experiments we use the relaxed Douglas-Rachford algorithm \eqref{RRR}. All variables are initialized to random real values between $0$ and $1$. In each example we first choose a value of the hyperparameter $\beta$ that works for that particular problem ($\beta=0.5$ for the dominating sets example, $\beta=0.8$ for the Boolean generative networks example), then choose a challenging instance of the problem and demonstrate the benefits of tuning the metric parameters while keeping the hyperparameter fixed. The first experiment will demonstrate how the tuning works and why it is important to have a small value of $\alpha$. The second experiment will demonstrate the main goal of tuning: preventing the algorithm from getting stuck on problems where such trapping is common for the untuned algorithm.

\subsubsection*{Experiments with dominating sets}
\begin{figure}[t]
  \centering
  \includegraphics[width=.4\textwidth]{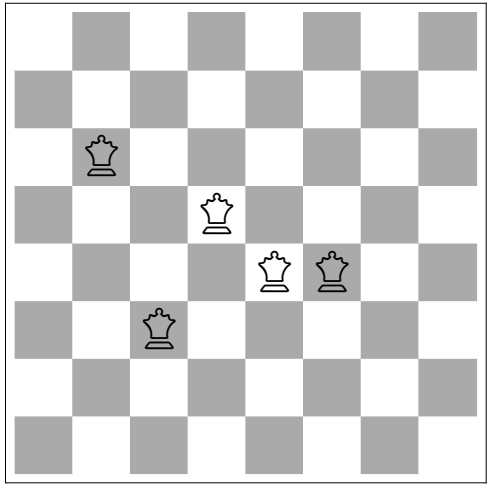}
  \caption{Five queens ``dominating" the $8\times 8$ chess board; that is, each unoccupied square is attacked by a queen. The domination number of the order 8 queens' graph is 5 because this is not possible with fewer queens~\cite{Cockayne1990}.}
  \label{fig:fivequeens}
\end{figure}

Consider a graph $G$ with vertices $V$ and (undirected) edges $E$. A subset $D\subset V$ \textit{dominates} $G$ if for every vertex $i\in V$ either $i\in D$ or there exists an adjacent vertex $j\in D$ such that $(i,j)\in E$. Finding dominating sets of minimum size $|D|$ is a well known NP-hard problem. Chess players will recognize the configuration shown in Figure \ref{fig:fivequeens} as an instance of a dominating set. Here the squares of the board are graph vertices and two squares are ``connected" by an edge if a queen placed on one attacks the other. The \textit{domination number} of this particular board/graph is $|D|= 5$ because the five queens attack all the other squares and this is not possible with fewer queens~\cite{Cockayne1990}.

That a vertex may be dominated either by itself or an adjacent vertex calls for two types of variable in a constraint formulation. A metric parameter should be introduced to control their relative ``weight." Our formulation uses vertex variables $y_i$, $i\in V$ and two variables for each edge $(i,j)\in E$ denoted $x_{i\to j}$ and $x_{j\to i}$ corresponding to the set of doubled (directed) edges $E_2$. One can think of $x_{i\to j}$ as a copy of vertex $i$ that vertex $j$ uses to express its domination status. Since there are only two variable types, a single metric parameter $\eta$ suffices to characterize the weight of the vertex variables relative to the edge variables. As explained in section \ref{sec:hypervsmetric}, we choose to absorb the metric weight by variable rescalings, specifically, in different discrete settings for the edge and vertex variables.

Constraint $A$ demands that every vertex is either ``dominating" or is ``dominated" by at least one adjacent vertex via an incident edge, with at most $|D|$ vertices being dominating:
\begin{equation*}
\begin{array}{lllr}
A:&\\
& \forall\; j\to i\in E_2: &\quad x_{j\to i}\in \{0,1\}\;,&(A1) \\[10pt]
& \forall\; i\in V: &\quad y_i\in \{0,\eta\}\quad \wedge\quad y_i/\eta\;+{\displaystyle \sum_{j\to i}}x_{j\to i}\;\ge 1\;,&\qquad (A2) \\[10pt]
& {\displaystyle \sum_{i\;\in\; V}} y_i/\eta \,\le\, |D|\;.&&(A3)\\[20pt]
\end{array}
\end{equation*}
Constraint $B$ demands that all edge variables agree with their associated vertex variable:
\begin{equation*}
\begin{array}{lllr}
B:& \forall\; j\to i\in E_2: &\quad x_{j\to i} = y_j/\eta\;.&
\end{array}
\end{equation*}
We project to the constraint sets using the metric $d((x_1,y_1),(x_2,y_2))=\sqrt{||x_1-x_2||^2+||y_1-y_2||^2}$. The following is pseudocode for projecting to the discrete set $A$:
\begin{itemize}
    \item For all $j\to i\in E_2$, set $P_A(x_{j\to i})=0$ if $x_{j\to i}<1/2$ and set it to $1$ otherwise. Also, set $P_A(y_i)=0$ for all $i\in V$. Call this the base state. This state satisfies constraint ($A$1).
    \item Determine which $j$ has the largest $x_{j\to i}$ for each $i\in V$. Call this $j_{max}(i)$. Also, for each $i\in V$ determine the change in squared distance $d_1(i)$, relative to the base state, for $i$ to be dominating, or $y_i=\eta$, as well as the change in squared distance $d_0(i)$, relative to the base state, for $i$ to be dominated, or $y_i=0$, $x_{j_{max}(i)\to i}=1$. Since only $y_i$ or $x_{j_{max}(i)\to i}$ might need to change, these are easy computations. The results are:
    \begin{align*}
    d_1(i)&=\eta^2-2\eta y_i,\\ d_0(i)&=\max(0,1-2x_{j_{max}(i)\to i})
    \end{align*}
    Setting $y_i=\eta$ whenever $d_1(i)<d_0(i)$ or $x_{j_{max}(i)\to i}=1$ whenever $d_1(i)>d_0(i)$ would be distance minimizing to constraint ($A$2). However, this may not satisfy ($A$3).
    \item To determine the distance minimizing set of vertices that should be dominating rather than dominated, we consider the differences $d_0(i)-d_1(i)$ for all $i\in V$. While vertices with $d_0(i)<d_1(i)$ should be dominated ($y_i=0$), of the vertices with $d_0(i)>d_1(i)$ at most $|D|$ may be dominating ($y_i=\eta$) by constraint ($A$3). If this number exceeds $|D|$, the distance minimizing subset of dominating vertices is those $|D|$ with the largest values of $d_0(i)-d_1(i)$.
\end{itemize}
Because of the vertex ranking in the last step, constraint $A$ is slightly non-local. Projection to set $B$ is much simpler:
\begin{itemize}
    \item For all $j\in V$, compute the weighted average $\left(\eta y_j + \sum_{j\to i} x_{j\to i}\right)/(\eta^2+\text{degree}(j))$, where $\text{degree}(j)$ is the number of edges leaving node $j$. Set $P_B(x_{j\to i})$ equal to this average, and set $P_B(y_j)$ to $\eta$ times this average.
\end{itemize}

\begin{table}
\begin{tabular}{ c|c|c } 
 $\alpha$ & successes/trials & iterations/solution\\
 \hline
 $0$ & $82/100$ & $4.18\times10^5$ \\ 
 $10^{-5}$ & $87/100$ & $3.59\times10^5$ \\ 
 $10^{-4}$ & $89/100$ & $3.45\times10^5$ \\ 
 $10^{-3}$ & $57/100$ & $1.18\times10^6$ \\ 
 \hline
\end{tabular}
\vspace{.1in}
\caption{Performance statistics for several values of $\alpha$ in searches for a size $6$ dominating set in the $12\times12$ queens' graph. Each trial was capped at $10^6$ iterations.}
\label{tab:ds}
\end{table}

Table \ref{tab:ds} gives performance statistics for several values of $\alpha$ in searches for a dominating set of size $|D|=6$ (the smallest possible~\cite{Cockayne1990}) for the queens' graph of order 12. The initial conditions were chosen randomly for each of the $100$ trials, but for the sake of fair comparison the $100$ sets of initial conditions were the same for each value of $\alpha$. In all cases we initialize $\eta=1$.

The main point of metric tuning is to avoid traps, not necessarily to speed up the algorithm for applications such as this in which trapping is rare; nonetheless, it is worth noting that the algorithm actually found solutions somewhat more quickly and reliably with $\alpha=10^{-5}$ and $\alpha=10^{-4}$ than with $\alpha=0$. However, more tuning is not necessarily better: $\alpha=10^{-3}$ does worse than $\alpha=0$.

\begin{figure}
  \centering
  \includegraphics[width=.97\textwidth]{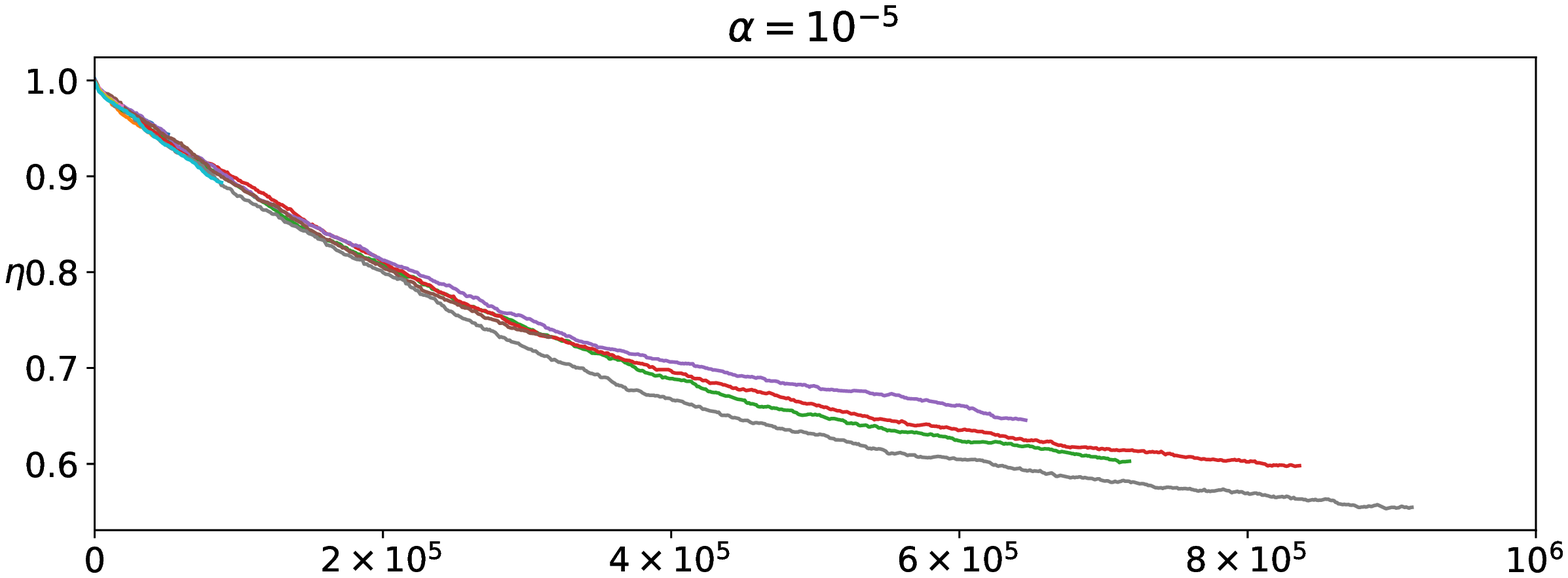}
  \includegraphics[width=.97\textwidth]{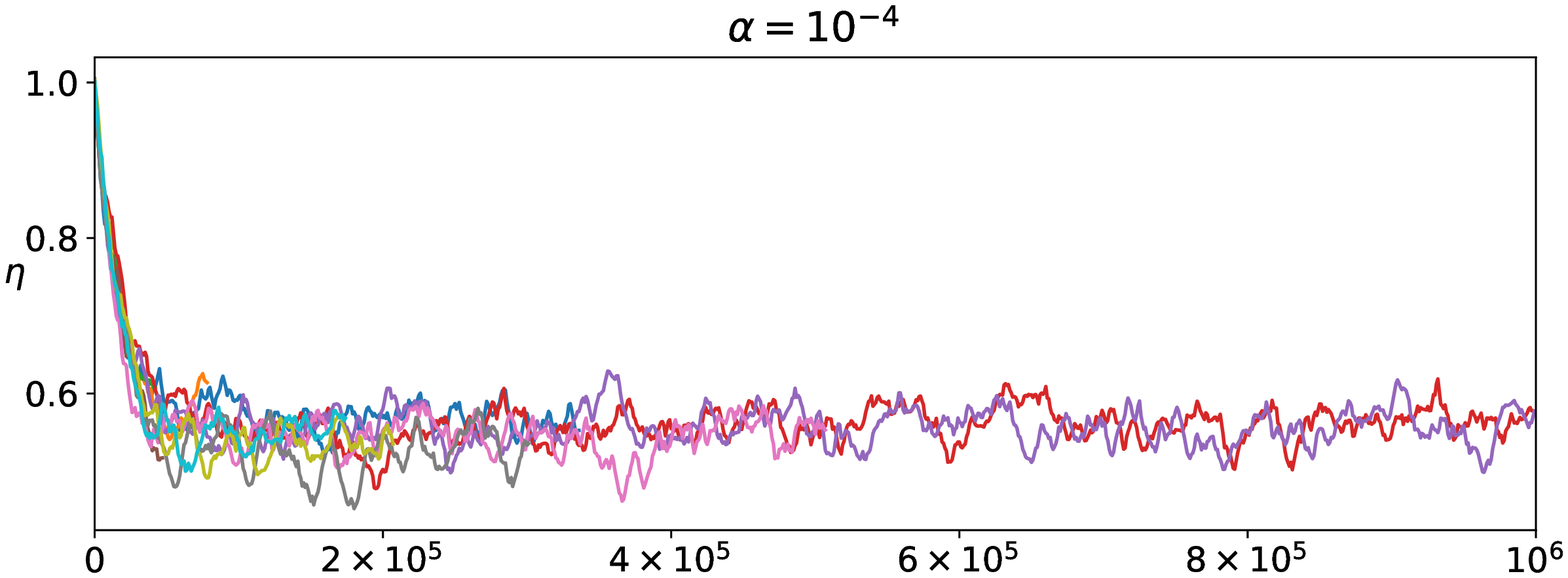}
  \includegraphics[width=\textwidth]{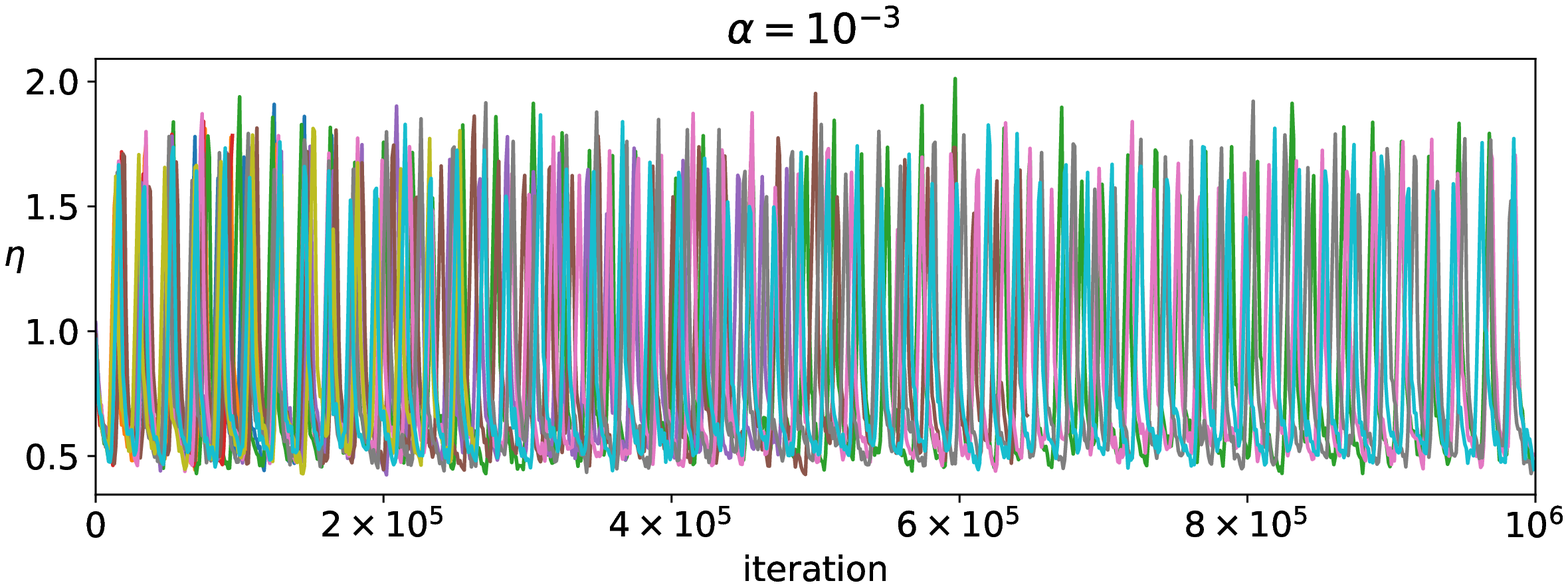}
  \caption{Time series of the metric parameter $\eta$ for ten trials at three difference values of $\alpha$, searching for a dominating set for the order 12 queens' graph. Each curve terminates when the algorithm finds a solution. }
  \label{fig:dstune}
\end{figure}

To understand the behavior with $\alpha$, it is helpful to plot the metric parameter $\eta$ as a function of iteration number. Figure \ref{fig:dstune} does this for ten trials of the three nonzero settings of $\alpha$. Each curve in the figure terminates when the algorithm finds a solution. For $\alpha=10^{-5}$, $\eta$ gently slides toward a value of around $0.6$, terminating early if a solution is found. With $\alpha=10^{-4}$ it reaches something near $0.6$ more quickly, though it tends to fluctuate more. A value somewhere between $\alpha=10^{-5}$ and $\alpha=10^{-4}$ might be optimal for this problem, but there is a generous range of $\alpha$ that works just fine: As long as $\alpha$ is significantly greater than $1/N$, where $N$ is the maximum number of iterations, there will be plenty of time for the tuning to take place. However, larger $\alpha$ is not always better, as we see with $\alpha=10^{-3}$. Here $\eta$ has large amplitude oscillations, almost periodically on a time scale of $10^4$ iterations. The metric is no longer evolving adiabatically and making a quasi-monotonic approach to a steady value. It is now changing so quickly that the local convergence properties of the Douglas-Rachford algorithm are compromised, resulting in the poor performance of $\alpha=10^{-3}$ in Table \ref{tab:ds}. This is why it is important to keep $\alpha$ small.

\subsubsection*{Experiments with Boolean generative networks}
In the previous experiment metric parameter tuning provided marginal benefits, slightly reducing the number of iterations needed to find a solution. In this next experiment we shall see an application in which trapping is common and the benefits of tuning are much more dramatic. In fact, tuning is essential to have any hope of finding solutions at all.

\begin{figure}
  \centering
  \includegraphics[width=.9\textwidth]{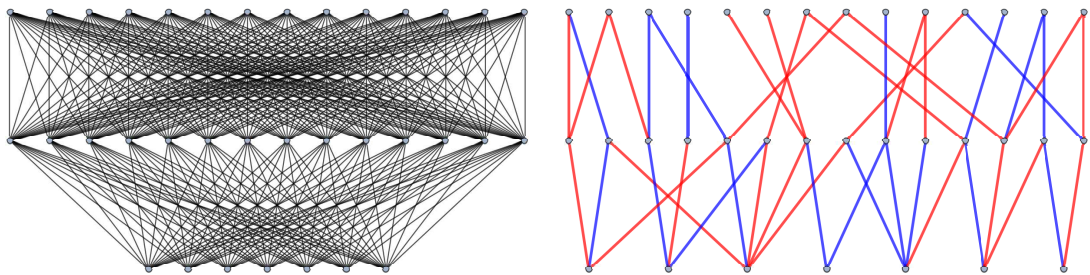}
  \caption{Network architecture (left) that takes 7 Boolean inputs and outputs 14 Boolean data; there is also an intermediate layer that holds 14 Boolean values. The trained BGN (right) is a logical circuit that uses relatively few of the available network edges as wires. All the nodes not in the input layer apply \textsc{OR} to the wires incident from below. Training also determines the placement of the \textsc{NOT} gates, which are indicated by red wires.}
  \label{fig:BGNcircuit}
\end{figure}

Our second example of automated metric parameter tuning is unsupervised machine learning with a Boolean generative network, or BGN \cite{BGN}. The data in this application consist of a set of $D$ Boolean strings of length $N$, and the network is tasked with discovering a Boolean circuit that generates the strings from a smaller number $M<N$ of Boolean ``latent variables" whose values for each data string are unknown. Figure \ref{fig:BGNcircuit} compares a network before and after training. Nodes of the network are arranged in layers and initially each node can potentially receive input from any node in the layer below. However, the data come with the promise or hypothesis that they can be generated with only \textsc{NOT} and 2-input \textsc{OR} gates. Through training the network must therefore discover which of the many edges are utilized in the circuit, that is, the ``wires" of the circuit, and whether the Boolean value is to be negated when traversing the wire.

The truth value at a node $i$ of the BGN is encoded by a variable $y_i$ at that node and also by copies of that truth value on all its out-edges, $x_{i\to j}$. These two variable types have the same interpretation they had in the dominating set application and serve to localize the constraints. However, by the uni-directional nature of the BGN logic, there is no need to have a second $x$ variable on each edge. Just as with dominating sets, the semantic equivalence of the $x$ and $y$ variables does not imply metrical equivalence and we control their relative scale with metric parameters $\theta$ for $x$ and $\eta$ for $y$.

\begin{figure}
  \centering
  \includegraphics[width=.5\textwidth]{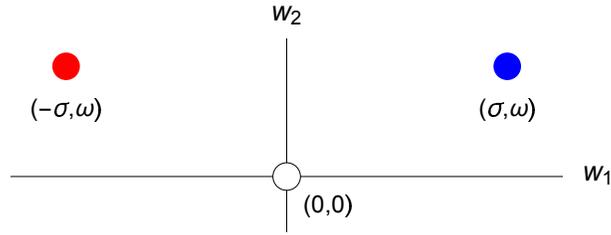}
  \caption{A 2D space $(w_1,w_2)$ is required to encode the three possible states of edges in a BGN. The red and blue points correspond to edges utilized as wires with and without negation; the open circle encodes the absence of a wire. Two metric parameters, $\sigma$ and $\omega$, define the isosceles geometry of the states.}
  \label{fig:wirestates}
\end{figure}

The relevance of metric parameters in BGNs is made especially obvious when we also consider the variables $w_{i\to j}$ used to encode whether an edge $i\to j$ of the network is a wire and if so, its negation state. Not only is the wire-status of an edge semantically distinct from the truth states it operates on, the existence of a wire (on the edge) should have independence, metrically, over its two negation states. To encode the three states of an edge (no wire, negating wire, non-negating wire) as three points we give each $w$ variable two independent components $w=(w_1,w_2)$ and represent the wire states with the vertices of an isosceles triangle as shown in Figure \ref{fig:wirestates}. This is the most general metrical relationship among the states that respects the symmetry between negated and non-negated wires. The metric parameter $\omega$ now controls the distance between presence and absence of a wire, while $\sigma$ controls the distance between the presence and absence of a \textsc{NOT}.

Here is a brief overview of the two constraint sets; see \cite{BGN} for further discussion of the projections in the BGN problem. In order to solve the problem we must create one copy of the network for each of the $D$ data strings. The $A$ constraint demands that every node, in each network copy, makes sense locally: 
\begin{itemize}
    \item Each incoming $w$ variable takes on one of the three values in Figure \ref{fig:wirestates} (representing a negating wire, non-negating wire, or no wire).
    \item Each incoming $x$ variable is either $0$ or $\theta$ (signifying that it carries value \textsc{False} or \textsc{True}).
    \item The $y$ variable at the node is $0$ if all of the incoming wires pass \textsc{False} or $\eta$ if at least one wire passes a \textsc{True}. 
\end{itemize}

The $B$ constraint makes the various copies agree:
\begin{itemize}
    \item The $x$ variable on an edge represents a copy of the truth value at its lower node in the network, and all the outgoing $x$'s from that node must agree with the $y$ variable at that node. These variables need not take discrete values in this equality constraint because the $A$ constraint sees to that. The distance-minimizing way to make the $x$'s and $y$'s agree is a weighted average, just like in the $B$ projection for the dominating sets problem.
    \item The $w$ variables represent the wire states, which, in a solution, must be the same across all $D$ copies of the network. Again, in this equality constraint the variables may take continuous values, so we simply average each edge's $w$ variable over all copies of the network.
\end{itemize} 

We will refer to tuning the four metric parameters $\omega$, $\sigma$, $\theta$, and $\eta$ as tuning by type. Now, it may happen that certain nodes and edges in the network are more important than others. For instance, nodes in different layers or having different in- or out-degrees could play markedly different roles. It is natural then to let the metric parameters be different for every node and edge in the network, promoting $\omega$, $\sigma$, $\theta$, and $\eta$ to $\omega_{i\to j}$, $\sigma_{i\to j}$, $\theta_{i\to j}$, and $\eta_i$. We will refer to this as tuning by type and location. One might also argue that the $D$ data strings have differing inherent difficulty and deserve suitably tuned metric parameters applied to their copies of the network. Instead of a single $\omega_{i\to j}$ parameter for edge $i\to j$ there would then be parameters $\omega_{k,\,i\to j}$ for $k=1,\ldots, D$ and similarly for the other variable types. We will refer to this as tuning by type, location, and data item.

We compared the different degrees of metric parameter tuning on a synthetic data set generated by the circuit in Figure \ref{fig:BGNcircuit} with $M=7$ inputs and $N=14$ outputs. For training we used all $D=87$ unique data strings generated by this circuit, initialized all metric parameters to $1$, and initialized the variables (on nodes and edges of the complete architecture) to random values between $0$ and $1$. Solutions (circuits) were only required to generate all $D$ data strings for some setting of the $M$ inputs. Typically the solution circuits additionally generated strings not among the $D$ data strings, though never a set of size $2^M$.

\begin{table}
\begin{tabular}{ l|c|c|c } 
 tuning & successes/trials & iterations/solution & iterations/second \\
 \hline
 none & $0/100$ & --- & $337.77$ \\ 
 by type & $2/100$ & $5.07\times10^6$ & $396.56$ \\ 
 by type and location & $58/100$ & $1.85\times10^5$ & $427.73$ \\ 
 by type, location, and data item & $56/100$ & $1.89\times10^5$ & $198.03$ \\ 
 \hline
\end{tabular}
\vspace{.1in}
\caption{Performance statistics for four degrees of metric parameter tuning in an instance of the BGN problem (Figure \ref{fig:BGNcircuit}). Each trial was capped at $10^5$ iterations.}
\label{tab:bgn}
\end{table}

Table \ref{tab:bgn} gives performance statistics for the four degrees of tuning, each run for $100$ trials capped at $10^5$ iterations per trial. As in the previous experiment, the random initial conditions were the same for each of the four tuning approaches, but different from one trial to the next. Clearly no tuning at all is not a viable option. Tuning by type is better, but still rarely succeeds within the limit of $10^5$ iterations per trial. Evidently it is necessary to tune by type and location, meaning that there is a $\theta_{i\to j}$ for each edge, an $\eta_i$ for each node, and so on. If tuning by type alone had been more effective, one might have considered doing so by hand: trying some choice of metric parameters for one trial, tweaking them according to which variable types appear to be getting stuck, and repeating. But tuning by type and location involves hundreds of metric parameters, making it infeasible to tune by hand. This is what makes our automated approach to tuning helpful.

The additional refinement of tuning by type, location, and data item does not make the search any more successful, and indeed it comes with a considerable slowdown in iterations performed per second. This is the drawback of tuning a large number of metric parameters. If the metric parameters are so numerous that tuning-work becomes an issue, one should first ask if having so many is truly necessary. In the present example, tuning by data item did not add any benefit that was not already seen with tuning by type and location. In an application in which the many metric parameters are all necessary one should consider ways of reducing the work, such as applying the metric updates less frequently than after every iteration.

\begin{figure}
  \centering
  \includegraphics[width=.9\textwidth]{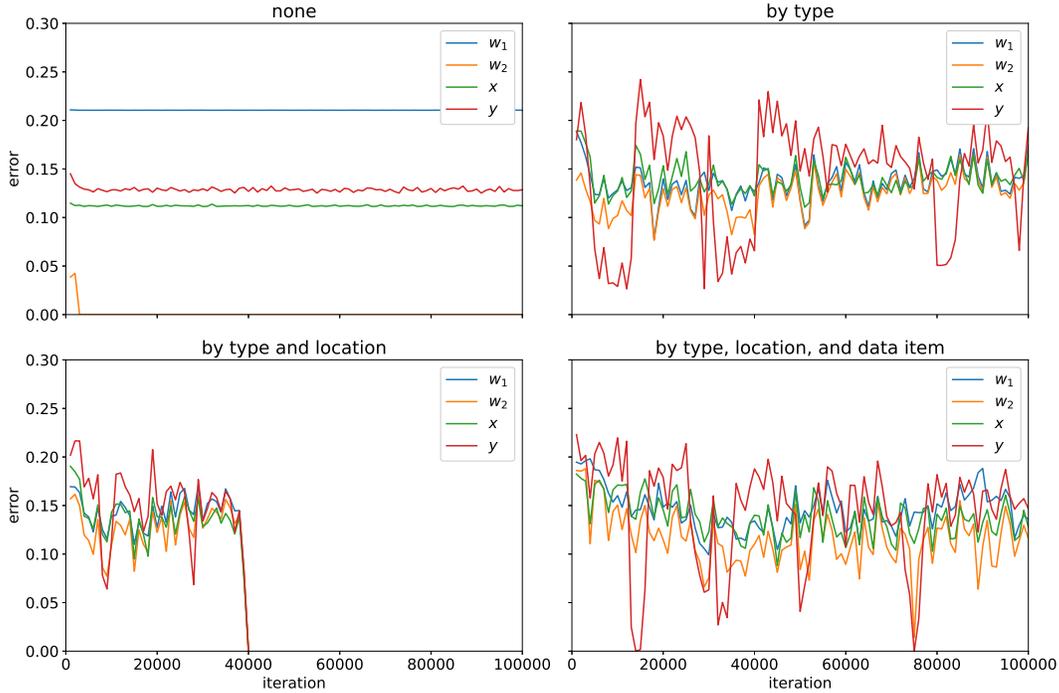}
  \caption{Time series of the constraint errors for a single run of the BGN problem with four levels of metric parameter tuning. \textit{Top left:} no tuning. \textit{Top right:} tuning by variable type. \textit{Bottom left:} tuning by type and location within the network. \textit{Bottom right:} tuning by type, location, and data item. The four errors in each plot are those of the four variable types described in the text.}
  \label{fig:bgn}
\end{figure}

For insight on how the four degrees of tuning affect the behavior of the Douglas-Rachford algorithm, it is helpful to look at plots of constraint error vs. iteration number. This is shown in Figure \ref{fig:bgn} where the error is broken down into contributions from each of the four variable types: $\epsilon_{w_1}$, $\epsilon_{w_2}$, $\epsilon_{x}$, $\epsilon_{y}$. The ideal behavior would be all four types fluctuating more or less in concert and with about the same amplitude, indicating that the difficulty of solving the problem is being shared equitably among the variable types.

What we see in the top left plot is far from the ideal behavior. This plot is for no tuning ($\alpha=0$) with the naive parameter choice $\sigma=\omega=\theta=\eta=1$. The constraint error for $w_2$ (wire vs. no-wire) quickly drops to zero, but the others stay at roughly constant nonzero values. Clearly this is an instance of the burden of constraint satisfaction not being shared among all the variable types. The result is that the search is stuck on a non-solution.

In the top right plot, we again initialize the metric parameters to $1$ but now we set $\alpha=10^{-4}$ to allow the four parameters to be slowly tuned as the search progresses. One can see that the search does not get stuck the way it did with $\alpha=0$, but the four variable types do not always move in concert. The $w_{1,2}$ and $x$ errors become highly correlated, but the $y$ error vacillates between much larger values and much smaller values while the others stay relatively stable. This could be in part because there are fewer $y$ variables than there are $w$'s or $x$'s, since $y$'s are associated with nodes whereas $w$'s and $x$'s are associated with edges. It is natural to expect the relative fluctuations to be larger for the smaller group of variables.

Nonetheless, we can still improve the correlation. In the bottom left plot, we promote $\omega$ to $\omega_{i\to j}$, $\sigma$ to $\sigma_{i\to j}$, and so on, to allow the metric parameters to vary also by location within the network. We still initialize the metric parameters to $1$ and use $\alpha=10^{-4}$. Here the four constraint errors appear more strongly correlated and the $y$ error no longer fluctuates with much greater amplitude than the others. The search terminates early because it finds a solution after about $4\times10^4$ iterations.

Finally, in the bottom right we go one step further, promoting $\omega_{i\to j}$ to $\omega_{k,\,i\to j}$, $\sigma_{i\to j}$ to $\sigma_{k,\,i\to j}$, and so on, to allow the metric parameters to vary by data item as well. As always we initialize the metric parameters to $1$, and we still use $\alpha=10^{-4}$. The correlation is certainly better than it was with no tuning and somewhat better than with tuning by type alone, but no better than tuning by type and location. 

We have also tried other tuning combinations; for instance, tuning by type and data item but not location performed no better than tuning by type alone. Among the many possible combinations, we conclude that tuning by type and location makes the greatest difference in the algorithm's performance.

One possible interpretation of these results is that metric parameters have tangible effects only when they target systematic or structural properties of the variables. The nodes of the BGN distinguish themselves by their in- and out-degrees, and also their distance from data constraints (imposed only at the output layer). We should therefore expect a metric-tuning benefit based on network location just as for the more obvious case of variables distinct by type (node vs. edge). In the case of variables differing only by the data index, the structural bias is much weaker. For most of these variables the structural difference is only indirect, through the fixed-value data constraints at the output nodes.

An alternative viewpoint is that metric parameters facilitate symmetry breaking behavior. For example, it may be important that  nodes within the same layer are able to develop unique characteristics, even while the associated node and edge variables have a permutation symmetry. Likewise, some data items might pose a greater challenge to the circuit than others, and the metric parameters of their associated variables would serve to break that form of permutation symmetry. Although our survey of results is far from comprehensive, the absence of a noticeable benefit from tuning by data index leads us to believe that the structural hypothesis is better supported than symmetry breaking.

\section{Conclusions}

We suspect that quite a few applications of iterative projection methods on nonconvex problems may have been abandoned because of a failure to recognize a sensitivity to parameters. Even in the case of hyperparameters, a setting in a range that is ``safe" with respect to local convergence (where the constraints may be approximated as convex) overlooks the fact that these parameters  also have a profound effect on the global characteristics of the search. We saw an example of this in the Douglas-Rachford generalization of section \ref{sec:doublereflect}, where $\delta$ needs to be tuned to a sweet spot for effective search. 

The tuning of metric parameters is equally important in applications where there is a rescaling freedom among variables not subject to symmetry. Such applications, without the translation or permutation symmetries of phase retrieval or a sudoku puzzle, are relatively new for iterative projection methods. Recognizing the role of the metric in such applications, and not arbitrarily assigning 1 as the relative scale, is an important first step. One can tune the metric parameters by hand, but here we have introduced a method for updating them automatically based on constraint errors. Actively updating the metric during the search does not disturb the convergence of the algorithm, as long as the updates are adiabatically slow. The automatic update approach is particularly useful in applications in which the naive metric fails and there are too many metric parameters to tune by hand.

\vskip 6mm
\noindent{\bf Acknowledgments}

\noindent The authors thank Avinash Mandaiya for noticing the two variable types in the dominating set problem and Jim Sethna for useful conversations.

\end{document}